\documentclass[12pt]{article}
\usepackage{epsfig}
\title{Conway's Nightmare: Brahmagupta and Butterflies}
\author{Richard Evan Schwartz \thanks{Supported by N.S.F. Grant DMS-2102802}}

\def\endproof{$\spadesuit$  \newline}

\def\R{\mbox{\boldmath{$R$}}}%

\begin{document}

\noindent
  \maketitle

  A {\it cyclic quad\/} is a convex quadrilateral
  whose vertices all lie on the same circle.
  Equivalently, opposite interior angles sum to $\pi$. That is, if we let
  $\alpha,\beta,\gamma,\delta$ be the successive
  interior angles, then
  $\alpha+\gamma=\beta+\delta=\pi$.

\begin{center}
\resizebox{!}{1in}{\includegraphics{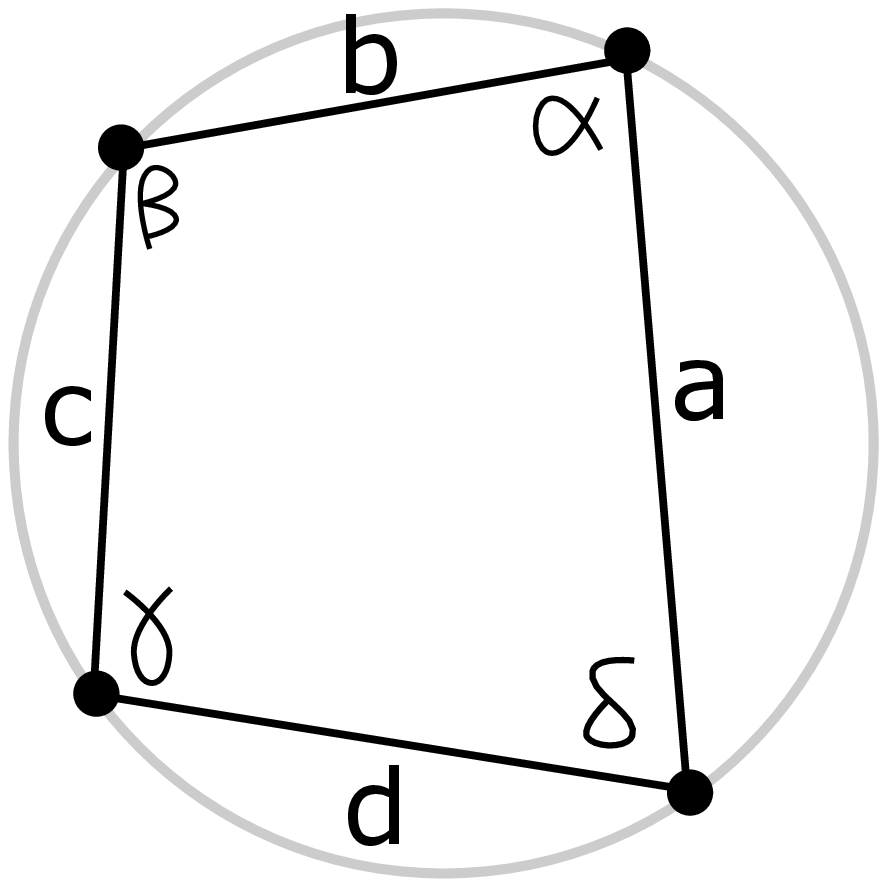}}
\end{center}

  Let $A$ be the area of such a quadrilateral and let
  $a,b,c,d$ be the side lengths.  Let \begin{equation}
  B=\sqrt{(s-a)(s-b)(s-c)(s-d)}, \hskip 30 pt s=\frac{a+b+c+d}{2}.
\end{equation}
Brahmagupta's formula says that $A=B$. I prefer to write
    \begin{equation}
      A^2=B^2=(s-a)(s-b)(s-c)(s-d).
    \end{equation}

Brahmagupta's formula goes back $1400$ years.
The historical paper \cite{KIC} discusses
how Brahmagupta himself may have
proved it.  John Conway long sought a simple and beautiful
{\it geometric\/} proof of Brahmagupta's formula,
like Sam Vandervelde's recent proof \cite{VAN}.
Is Brahmagupta's formula really a geometric
result?  I am not so sure.   When generalized
to polygons with more sides, as in
\cite{ROB}, the discussion turns
decidedly and deeply algebraic.   In any case,
I am less interested in a geometric proof than I am
in a proof that is short, conceptual, and without calculation.

My proof only depends on basic facts about polynomials
and continuity, but I got the idea thinking about things in modern
mathematics like flat cone surfaces, control theory, and ergodicity.
(It also helps to have smart friends; see the
acknowledgements at the end.)
I call this proof ``Conway's nightmare''
because Sam had been calling his proof
``Conway's dream'' in an early draft
of his paper.
My proof probably would not have satisfied
Conway.
I will set the stage for the proof, give the proof, then discuss
the mathematics that inspired it.
\newline
\newline
{\bf Setting the Stage:\/} Let me
explain my favorite proof of
the Pythagorean Theorem.   If we {\it fix the
angles\/} of a right triangle $T_c$ with sides $a,b,c$
and hypotenuse $c$, then ${\rm area\/}(T_c)$ is proportional
to $c^2$ because $a$ and $b$ are both
proportional to $c$.
We might write this as ${\rm area\/}(T_c) \propto c^2$.
\begin{center}
\resizebox{!}{.6in}{\includegraphics{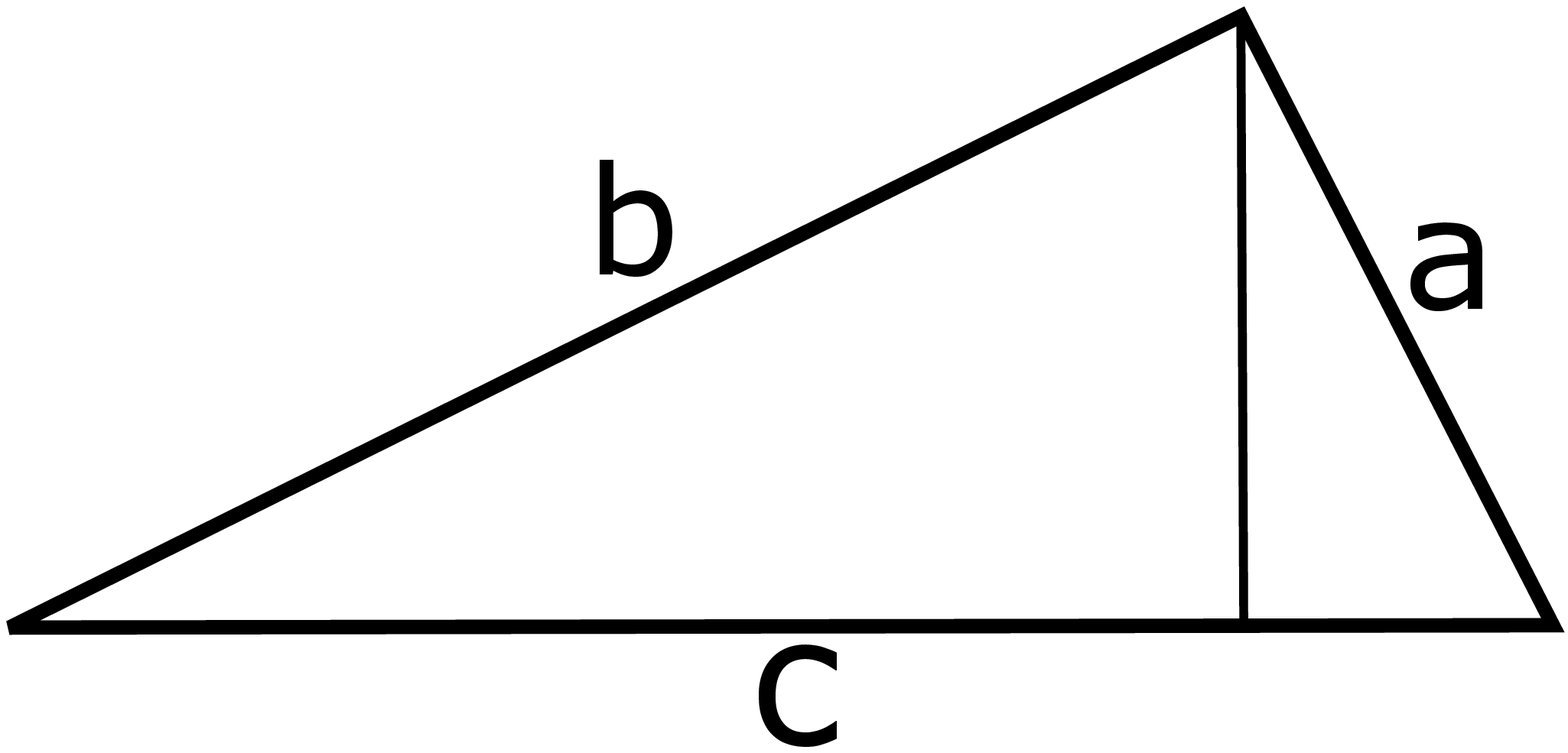}}
\end{center}
Put another way, ${\rm area\/}(T_c)=\lambda c^2$, where
$\lambda$ is a constant that depends
on the angles.   We can divide $T_c$ into two
smaller triangles $T_a$ and $T_b$ by
dropping the altitude perpendicular to $c$.
Each smaller triangle has the same angles as $T_c$, so
the same $\lambda$ works for all
$3$ triangles.  Since
${\rm area\/}(T_c)={\rm area\/}(T_a)+{\rm area\/}(T_b)$ we have
$\lambda c^2=\lambda a^2 + \lambda b^2$.
Cancelling $\lambda$ gives the result.

In this proof, we considered
what happens when we vary a right triangle in a special
way: changing its size/position without changing the angles.
Let us call this operation {\it morphing\/}.  The name will
be more apt when we apply it to cyclic quads.
Noting how the relevant quantities associated to a right triangle
change when we morph, we 
get the desired relation up to a constant that cancels out.  This
proves the Pythagorean Theorem one ``angle-type'' (a.k.a. {\it
  similarity class\/})
at a time.

My proof of Brahmagupta's formula has the same flavor.
There is one part (second paragraph) that examines
how the relevant quantities vary when we morph. -- i.e., vary
without changing the angles.
This analysis alone gives a useful partial result.
The other part of the proof (first paragraph) combines the
morphing result with an obvious result concerning
another operation, {\it recutting\/}, to
close the deal.

Let me recall two notions that arise in my proof.
The {\it signed distance\/}
between two points $p,q$ on the same line $\ell$ is
the dot product $(p-q) \cdot u$, where $u$ is one of the two unit
vectors
parallel to $\ell$.  The {\it signed area\/} of a polygon
is the sum of the areas of all the regions it bounds, weighted
by the number of times the polygon winds around
the points in each region.  This quantity also has an
algebraic expression, in terms of determinants.  Both
quantities require choices; this amounts to choosing signs
for some example.
\newline
\newline
\noindent
{\bf The Proof:\/}
Let $C=A^2/B^2$.  Let $X$ be the space of
cyclic quads.
To {\it morph\/} a quad is to
replace it by one with the same angles.  To {\it recut\/} a
quad is to cut along a diagonal and reverse one
triangle.  Recutting preserves $C$. 
We claim morphing does too.  From any point in $X$
we can reach all nearby points by morphs and recuts.
(Recut, morph, re-recut to perturb one pair of
opposite angles; repeat using the other diagonal;
morph one final time.)
Hence $C$ is constant on $X$.  Since $X$ has squares, $C=1$.

\begin{center}
\resizebox{!}{1.9in}{\includegraphics{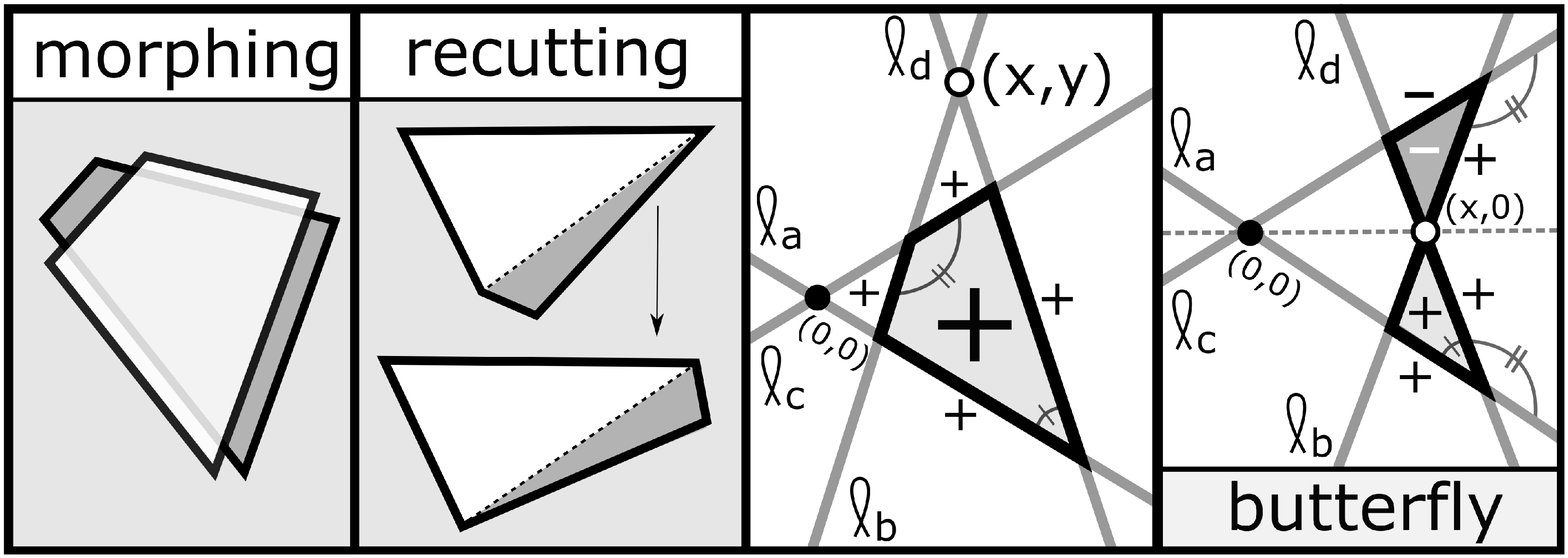}}
\end{center}

\noindent
\underline{Proof of Claim:}
Let $L(\rho,\sigma)$ be the space 
of lines $\ell_a,\ell_b,\ell_c,\ell_d$ with slopes
$\rho,\sigma,-\rho,-\sigma$ and
$\ell_a \cap \ell_c=(0,0)$.  Parametrize $L \cong \R^2$ by
$(x,y)=\ell_b \cap \ell_d$.
Let $a,b,c,d$ be the
signed distances between
vertices of the associated quads, and let $A$ be the signed
area.  Choose the signs
so that $a,b,c,d,A>0$ in a convex case.
$B^2(x,y)$ is a degree $4$ polynomial as
$a(x,y),..., d(x,y)$ are linear. $A(x,y)$, a
sum of determinants of linear functions,
is a degree $2$ polynomial.
If $xy=0$ the  quads are {\it butterflies\/},
so $A=0$; also $|a|=|c|$ and $|b|=|d|$ and $ab=-cd$, so
$2$ factors of $B^2$ vanish.  Given their
degrees, $A(x,y) \propto xy$ and
$B^2(x,y) \propto (xy)^2$.  Hence
$C|_L$ is constant.  Our claim follows: A (generic) quad and its
morphs are all isometric to quads in the same $L$.
\endproof

\noindent
    {\bf Discussion:\/}
    The main idea is that the
    roots, counted with multiplicity, determine a
    real polynomial up to constants provided
    that the number of roots equals the degree of
    the polynomial.
    A single evaluation then determines the constant.
    My proof can be summarized like this: The
    function $C=A^2/B^2$ is invariant under
    recutting.  It is also invariant under morphing because,
    when analytically continued, $A^2$ and $B^2$
    vanish to the same order on the set of
    butterflies and nowhere else.  The
    recutting/morphing process spreads the constancy of $C$
    through $X$ like a virus.

    To make this idea work, we have to enlarge the space $X$ so that
    it includes some nonconvex quads, especially butterflies.
    We'll explain it from another point of view here.
   First of all, let us modify $X$ so that we
consider cyclic quads modulo isometry.
    Call quads {\it cousins\/} if they are morphs of each other.
We think of $X$ as a fiber bundle,
where the fibers are the cousin families.
Each fiber is a convex cone in $\R^2$.
We create a new space $X^*$ by
replacing these cones by the copies of
$\R^2$ which contain them.
The space $X^*$ is a
plane bundle with the same base.
The fibers are our $L$ spaces.

The fact that $A$ is a
degree $2$ polynomial on the fibers is
a key idea of Bill Thurston's
paper {\it Shapes of Polyhedra\/} \cite{THU}.
In Thurston's work, he introduces
local complex linear coordinates on the space
of flat cone spheres with prescribed
cone angles. Prescribing the cone angles
is like restricting to a fiber.  Thurston's coordinates
are like my $(x,y)$ coordinates. He shows
the area of a flat cone surface with
fixed cone angles is the diagonal
part of a Hermitian form in his coordinates.  There is also
a real valued version of this theory which
is even closer to my proof, exposited
recently in the A.M.S. {\it Notices\/} \cite{CAL} by
Danny Calegari.  The same ideas also arise in translation surfaces.
The main  point is that if you fix the slopes
of the lines (or the cone angles, in
Thurston's case), various algebraic functions
are simplified and become linear.

Now we discard $X^*$ and go back to $X$.
Once we know that $C$ is fiberwise constant
on $X$ how do we get $C=1$?  Let me mention
two alternate approaches first.
Each fiber contains triangles (i.e. degenerate
  quads) and then we could deduce $C=1$ by Heron's formula,
  a degenerate version of
  Brahmagupta's formula which is somewhat easier to prove.
  (There are other reductions of Brahmagupta to Heron,
  e.g. \cite{AH}.)  However, you would
  still need to prove Heron's formula.
Better yet, Peter Doyle noticed that each
fiber contains a (perhaps nonconvex) quad whose diagonals are
perpendicular, and that for such quads Brahmagupta's
formula can be verified with some clever but ultimately
easy algebra.  I'll leave this as a challenge.

My inspiration for the morphing/recutting proof came from
control-theory flavored proofs of ergodicity.
The prototypical example is
E. Hopf's proof \cite{HOP} that
the geodesic flow on a
hyperbolic surface is {\it ergodic\/}, meaning
that any invariant (measurable) function is
(almost everywhere) constant.
This may seem far-fetched, but consider the picture.

The geodesic flow lives
on a $3$-manifold, the unit tangent bundle of
the surface.  This $3$ manifold has $2$ invariant
codimension $1$ foliations,
the {\it stable foliation\/} and the
{\it unstable foliation\/}.    
The first step in Hopf's proof is to
use the expansion/contraction properties of
the flow along the leaves of the foliations
to establish the (almost everywhere) constancy on each leaf
of the stable foliation and each leaf of the unstable foliation.
The next step is to walk around,
going from a stable leaf to an unstable
leaf to a stable leaf, etc, to spread this
constancy around over the whole $3$-manifold.
Our space $X$ has $2$ codimension $1$ foliations,
each consisting of
quads sharing a pair of opposite interior angles.
The method I gave, in the first paragraph of the proof,
for moving around in $X$ is the same kind of
alternating walk through the leaves of these foliations.

I can't resist explaining $2$ other ways I might have done
the control theory part.   The first
approach minimizes ``effort'' and the second
approach minimizes exposition length, but
both make extra demands on the reader.

\begin{enumerate}
\item Using a finite sequence of morphs and recuts one
  can start with some point in $X$ and reach an open subset
  of $X$.  Hence $C$ is constant on an open subset of $X$.
  Since $X$ is connected and $C$ is an analytic function,
  $C$ is constant on $X$.  Since $X$ has squares, $C=1$.

\item One can join an arbitrary cyclic quad to a square
  by a finite sequence of morphs and recuts. Hence $C=1$.
\end{enumerate}

Here is an algorithm for the second approach.
To make it as clean as
possible we note that both morphing and recutting extend to
degenerate cyclic quads -- i.e.
triangles with one marked point.
{\it Morph so as to maximize the intersection angle
between the diagonals, recut along the longest
diagonal, repeat until done.\/}
Let $\phi=(d_1^2+d_2^2)/D^2$ where $d_1,d_2$ are
the diagonal lengths and $D$ is the diameter of the circle
containing the vertices.  The algorithm
produces a finite number of marked triangles, increasing
$\phi$ by a factor of at least $3/2$ each step, until
$\phi>1$.  Then we get one or two quads with perpendicular
diagonals, the last being a square.
\newline
\newline
{\bf Acknowledgements:\/}
Peter Doyle rekindled my interest in Brahmagupta's formula
by showing me Sam Vandervelde's proof.   Then Peter went
on to explain how one can rotate a cyclic quad so
that it is a quad of $L(\rho,\sigma)$.  Something
about this rang a bell, and once I realized that this was
just like Thurston's paper the rest fell into place.
Peter's insight about the connection
to $L(\rho,\sigma)$ was my main inspiration.
After I explained
my proof to Jeremy Kahn, he
suggested the idea of taking the quotient by translation
and working in $\R^2$ rather than $\R^4$.  This simplified
the algebra.  Danny Calegari, Dan Margalit,  Javi Gomez-Serrano, and
Joe Silverman all had helpful expository suggestions.

\end{document}